\documentclass[a4paper,12pt,fleqn]{article}

\usepackage{amssymb}
\usepackage{amsmath}
\usepackage{graphicx}
\usepackage{color}
\usepackage{eurosym}
\usepackage{txfonts}
\usepackage{courier}
\usepackage{pifont}

\usepackage[superscript]{cite}

\usepackage[colorlinks=false,urlcolor=blue,ps2pdf=true,
pdfpagemode=UseNone, pdfauthor={Robin Whitty},
pdfsubject={Mathematics of OR}, pdfcreator={Robin Whitty},
pdfproducer={Robin Whitty}]{hyperref}

\begin{document}

\mbox{\ }

\vspace{.6in}

\begin{center}
{\LARGE Some Comments on Multiple Discovery in Mathematics\footnote{Read at the 2014 Christmas Meeting of the British Society for the History of Mathematics}}

\bigskip
Version: \today

\vspace{.5in}
{\large Robin Whitty, Queen Mary University of London}
\end{center}

\vspace{.2in}
\noindent
\begin{quote}
{\small {\bf Abstract} Among perhaps many things common to Kuratowski's Theorem in graph theory, Reidemeister's Theorem in topology, and Cook's Theorem in theoretical computer science is this: all belong to the phenomenon of simultaneous discovery in mathematics. We are interested to know whether this phenomenon, and its close cousin repeated discovery, give rise to meaningful questions regarding causes, trends, categories, etc. With this in view we unearth many more examples, find some tenuous connections and draw some tentative conclusions.
}
\end{quote}

\vspace{.15in}
\noindent
The characterisation by forbidden minors of graph planarity; the Reidemiester moves in knot theory; and the NP-completeness of satisfiability, are all discoveries of twentieth century mathematics which occurred twice, more or less simultaneously, on different sides of the Atlantic.
One may enumerate such coincidences to one's heart's content. Michael Deakin wrote about some more\cite{Deakin} in his admirable magazine {\em Function}; Frank Harary reminisces\cite{Harary} about a dozen or so in graph theory, and there are whole books on individual instances, such as  Hall\cite{Hall} on Newton--Leibniz and Roquette\cite{Roquette} on the Main Theorem of modern algebra.

Mathematics is quintessentially a selfless joint endeavour, in which the pursuit of knowledge is its own reward. This is not a polite fiction; G.H. Hardy's somewhat sour remark\cite{Hardy} ``if a mathematician ... were to tell me that the
driving force in his work had been the desire to benefit humanity, then I should not believe him (nor should I think the better of him if I did)" nevertheless allows him to rank intellectual curiosity before ambition; and this ranking is engrained in a mathematician's upbringing. So simultaneity of mathematical discoveries tends to be tactfully overlooked: disputes about priority or naming or respective merits are in poor taste. The infamous affair of Erd\H{o}s and Selberg's simultaneous elementary proofs of the Prime Number Theorem, for instance, although meticulously documented\cite{Goldfeld}, is generally considered best forgotten. The topologist Joan Birman represented the profession eloquently\cite{Birman}  in relation to another fuss,\cite{Nasar} over who proved Poincar\'e's Conjecture: ``Reading [Nasar's biography of John Nash], I was proud of our decency as a community."

To take a less controversial example, Joseph Wedderburn's Little Theorem was proved first, marginally, and proved better by Leonard Dickson,\cite{Parshall} but there is no move to re-brand it Dickson--Wedderburn. Again, it is quite conceivable that Cauchy was aware of Galois' marginally prior discovery that prime order groups contain subgroups of the same prime order, it remains Cauchy's theorem for all that. But there is pause for thought here: Meo\cite{Meo} addresses the possibility that Cauchy indeed knew of but `disqualified' Galois's ideas on political or religious grounds. The historian of mathematics cannot quite overlook this: curiosity must be roused, or perhaps heckles raised.

And what are we to make of particularly spectacular examples of simultaneity? How could anyone shrug off the appearance, at the same moment and for the same purpose, but in utterly different cultures, of the Bernoulli numbers? They were essentially present in Faulhaber's Formula of 1631, which formula has a long prehistory in Greece, India and the Muslim world.\cite{Beery} But they became explicit in work of Jacob Bernoulli, posthumously published in Eidgenossenschaft Switzerland in 1713 and of Takakazu Seki, posthumously published in shogun Japan in 1712. Wimmer-Zagier and Zagier\cite{Zagier} examine the possible interactions between such widely separated mathematical cultures and propose that this is a genuine example of independence, and all the more striking because: ``Bernoulli numbers constitute a far more specific discovery [than determinants] that is not an inevitable consequence of any particular `general trend in the science.'''

Repeated discovery is by no means disconnected from its simultaneous cousin. Historians have declared (c.f. chapter~2 of Mancosu\cite{Mancosu}) that the Pappus--Guldin Theorems, relating displacement of centre of mass to surface areas and volumes of revolution, were independent discoveries by Pappus, around 300AD, and by early seventeenth century scholars, notably Paul Guldin. At least an analogy is suggested between the cultural divide separating exact contemporaries in Europe and Japan, and that separating fellow Europeans thirteen hundred years distant from each other.

An interesting instance of rediscovery is the classification of semiregular (aka Archimedean) tilings of the plane: this appears in chapter 2 of Kepler's 1619 {\em Harmonices Mundi} (in which also appears his third law of planetary motion) but there have been at least three rediscoveries, the most recent being simultaneous: in 1905 by the Scottish mathematician Duncan Sommerville (a near-neighbour and contemporary of Wedderburn)  and by the Italian entomologist Alfredo Andreini in 1907. In between there is a French publication dated 1887 by a Paul Robin whose identity seems obscure, although there was a chemist of that name active at the time who might legitimately have been interested in tilings. These rediscoveries might suggest multiple drivers for geometrical discovery: mathematics, philosophy (in the sense of natural philosophy) biology, and (maybe) chemistry.

So questions are, it would seem, being begged. But they are not, on the evidence, being answered in a systematic way. We might not expect anything as impressive as a Poisson model of discovery arrival, as has been proposed, tongue-in-cheek, by Ken Regan\cite{Lipton2} but we could hope for some kind of frame of reference or some parameters for discussion. Brian Blank, in his sweeping (and scathing) review\cite{Blank} of another Newton--Leibniz book, says  ``Historians and sociologists of science have long been fascinated with {\em multiple discoveries} ... Such discoveries are even more noteworthy when they exemplify the phenomenon of {\em convergence}---the intersection of research trajectories that have different initial directions." What is this `phenomenon of {\em convergence}' (Blank's emphasis)? Might it suggest a systematic framework in which to compare Newton (astronomicals) vs Leibniz (infinitesimals) with the Kepler--Robin--Andreini story? Such a framework would be an important and useful beast, but it appears to be an elusive one. Perhaps it is at work in Kuhn's famous examination\cite{Kuhn} of the history of energy conversation?
\begin{quote}``What we see ... is not really the simultaneous discovery of energy conservation. Rather it is the rapid and often disorderly emergence of the experimental and conceptual elements from which that theory was shortly to be compounded."
\end{quote}
`Rapid and often disorderly' might very well be said to apply to the emergence of the calculus. But in the mathematical sciences we might perhaps identify and describe and account for discoveries at a  finer level of granularity, and consequently
  more precisely and tidily, by focusing on individual milestones in mathematical progress, theorems and lemmas and their proofs, or conjectures which eventually acquire proofs.

There is a danger of drawing a false distinction: surely a theorem should not be isolated from the whole body of theory to which it belongs? Stedall\cite{Stedall} goes further: the history of mathematics should concern itself with
``... how, where, and why mathematics has been practised by people whose names will never appear in standard histories."\cite{Stedall}
In this view, a body of mathematical theory itself cannot be isolated from the times, cultures and processes which engendered it. Even the phrase `historians and sociologists of science', is to distinguish falsely. This was also the opinion of Robert Merton, pioneer of the sociology of science:
 \begin{quote}
 ``...for more than three centuries, there has been an intermittent mock battle between the advocates of the heroic theory and the theory of the social determination of discovery in science. In this conflict truth has often been the major casualty. For want of an alternative theory, we have been condemned to repeat the false disjunction between the heroic theory centered on men of genius and the sociological theory centered on the social determination of scientific discovery"\cite{Merton1}
 \end{quote}
 The quotation is from Merton's best known paper on multiple discoveries in science; the first ten pages of which are concerned with Francis Bacon's much earlier contemplations on the same subject, hence the reference to `three centuries'. To Bacon, Merton attributes the (implicit) identification of the causes of what he calls `multiples': ``the incremental accumulation of knowledge, the sustained social interaction between men of science  and the methodical use of procedures of inquiry." To which he adds Bacon's `Births of Time' thesis:
``Once the needed antecedent conditions obtain, discoveries are offshoots of their time, rather than turning up altogether at random." The Baconian view of simultaneity continues to be the prevailing one. Pickover, for instance, in his popular history\cite{Pickover} might find Georges Ifrah in disagreement:
``The true explanation lies in ... the profound unity of culture: the intelligence of {\em Homo sapiens} is universal and its potential is remarkably uniform in all parts of the world."
But this (and a literary flourish about numbers offering universal solace in times of turmoil) appear secondary to what must be taken to be his own view:
``Most likely, such simultaneous discoveries have occurred because the time was ripe for such discoveries, given humanity's accumulated knowledge at the time the discoveries were made."

Merton's investigations drew him to propound what has been called `Merton's Hypothesis':
``the apparently incorrigible hypothesis that singletons, rather than multiples, are the exception requiring distinctive explanation and that discoveries in science are, in principle, potential multiples." But the `distinctive explanation' is intended, of course, to be a sociological one, and Merton's sociology of scientific discovery is most interested in `reward systems' (indeed, a reprint of his paper appears in his collection\cite{Merton2} {\em The Sociology of Science} within a section called ``The reward system of science"). Subsequent work, notably Zuckerman's celebrated analysis of Nobel prizewinners,\cite{Zuckerman} has given much substance to this system. But while it might serve to illuminate the Newton--Leibniz spat (or indeed the Bernoulli--Bernoulli one which arguably prevented Jacob Bernoulli's work appearing somewhat before Takakazu Seki's) it does not appear to offer anything very tangible to the historically inclined mathematician, who will be more interested in how mathematicians discover, rather than how they behave. The historically inclined mathematician may after all, {\em pace} Merton and Stedall, see what can be done at the level of individual theorems, discovered by clever individuals. Neumann's examination\cite{Neumann} of the dating of Cauchy's Theorem, and Stigler's detective work\cite{Stigler} on Bayes' Theorem are  two compelling examples.

Our opening transatlantic simultaneity, Kuratowski's Theorem, has likewise been comprehensively examined by Kennedy et al.\cite{Kennedy} but rather in the matter of adjudicating on a priority dispute. It is one of the simultaneities cited by Harary\cite{Harary}, who endorses Merton's Hypothesis: ``independent discovery is probably more the rule than the exception, particularly in a field which is growing as rapidly as graph theory." A non-sociological `distinctive explanation' is certainly offered here: a new and rapidly expanding area in mathematics is likely to attract many clever new researchers, who are likely to uncover for themselves, independently, a certain subset of basic truths. But even more distinctive would be an account which took into account what was happening mathematically in Europe and in the US in the late 1920s. As Chartrand et al observe,\cite{Chartrand} Kuratowski's result concerned ``le probl\`eme des courbes gauches en topologies", while that of its US co-discoverers, Frink and Smith, self-effacingly reduced to a one-line abstract in {\em Bull. AMS}, concerned  ``Irreducible non-planar graphs". Orrin Frink and Paul Althaus Smith both completed their doctorates in 1926, in modern algebra at Columbia and in topology at Princeton respectively, but Frink published at the same time a new proof of a classic 1891 theorem of Julius Petersen on factorisation of regular graphs of odd degree. Kuratowski's theorem was also known to L.S. Pontryagin and his teacher P.S. Alexandrov in Moscow, and to Karl Menger in Vienna, all preoccupied with curves, surfaces, and the like, although Menger, like Frink, came to the problem of planarity with a graph-theoretic perspective, having become interested in map colouring and cubic graphs. This certainly seems to be one of Bacon's births of time: Jordan's theorem on plane curves dates from the 1880s but had been given what was at the time considered\cite{Hales1} to be its first complete treatment by Veblen in 1905, a point on a trajectory leading inevitably to planarity of ensembles of plane curves. But this is not meant to be a conclusive assertion but rather a provocative one!

Karl Menger makes an interesting point of contact with our second opening simultaneity: his study of cubic graphs dates from his membership of ``Reidemeister's Vienna seminar".\cite{Kennedy} Kurt Reidemeister was in Vienna between 1923 and 1925 and it was there that he was introduced to knot theory by Wilhelm Wirtinger. It seems plausible that a consideration of knot invariants (the fundamental group, for example) and questions about equivalence would combine with the Vienna seminar's coincidental study of graphs and maps to supply the necessary impetus for Reidemeister to discover his diagrammatic `moves'. For that matter, both the systematic enumeration of knots and the study of cubic graphs had been pioneered by the same man, Peter Guthrie Tait, in the 1880s. The simultaneous discovery of Reidemeister's theorem by James Waddell Alexander and his Princeton student Garland Baird Briggs is contrasted epistemologically with Reidemeister's by Moritz Epple;\cite{Epple} and seeing their respective diagrams side by side reinforces the idea that Reidemeister was influenced by graph theory, while Alexander's was a more purely algebraical approach. Certainly Alexander belonged to a `tradition' of topology, founded by Veblen and continued by Lefschetz.

Graph theory is a common thread taking us through the second war, with its enormous impact on discrete optimisation, to the 1970s by which time graph algorithms had become a driving force in the study of computational complexity.
The NP-completeness of satisfiability,  our final opening simultaneity, was proved in 1971 by Stephen Cook, then at the University of Toronto. This was certainly work which was `of its time': by 1971 Jack Edmonds was barely seventy miles distant in Waterloo and would have or could have discovered the same result; the same might be said of Richard Karp of whom Cook was a very recent Berkeley colleague. Cook moved to Toronto after falling foul of Berkeley's `reward system', so there is certainly sociological material to hand. Once again, however, there is a pressing mathematical claim on our attention: something roughly equivalent\cite{Lipton1} to Cook's Theorem was discovered simultaneously by Leonid Levin, a recent graduate and student of Kolmogorov, working at the Moscow Institute of Information Transmission. A fascinating confluence of mathematical themes is suggested: combinatorial optimisation developing in post-war USA and the Soviet Union, wonderfully invoked by Schrijver\cite{Schrijver} (a companion paper,\cite{Polyak} by Polyak, in the same journal issue gives a Soviet viewpoint); the theory of algorithms developing at Harvard (Michael Rabin, Alan Cobham) where Cook studied for his PhD; the Kolmogorov route to complexity (itself a `multiple', having been pre-invented by Ray Solomonoff, working in the US military research sector in the early 1960s, (c.f. p. 89 of Li and Vit\'anyi's book\cite{Li}), not to mention the evolution of computers, computer programming, the formal theory of languages and  artificial intelligence.

It is missing an opportunity, then, to explain away the simultaneity of Cook's Theorem  as a case of `the time was ripe'. Some simultaneities may have less milage in them, however. The discovery of the Max-Flow-Min-Cut Theorem of combinatorial optimisation, for example, is another transatlantic event, reported in a RAND report by Ford and Fulkerson in 1954 and independently, in a restricted version, by Anton Kotzig in his doctoral thesis of 1956, in Bratislava. But this is a multiple discovery which is quite unsurprising given mathematical and industrial preoccupations at the time. Indeed, the theorem is also  attributed to Claude Shannon and his colleagues who wrote in 1956 ``This theorem may appear almost obvious on physical grounds and appears to have been accepted without proof for some time by workers in communication theory."\cite{Aardal} It belongs perhaps to that most vague of categories of multiple discovery, the `folk theorem'!

Wedderburn's Little Theorem   must in hindsight appear to be equally inevitable, as a step in the progress of modern algebra. Moreover its co-discoverers and their colleagues were in close contact with each other. The main protagonists were all even at the same institution during the critical period. Parshall's analysis\cite{Parshall} is surely definitive, although the theorem itself has had a rich afterlife, acquiring alternative proofs and moonlighting in projective geometry. The two 1896 simultaneous (non-elementary) proofs of the Prime Number Theorem were inevitable after 1859 when Riemann set out the game-plan.\cite{Fine} Around 1980, Georgy Egorychev and Dmitry Falikman simultaneously proved the 1926 van der Waerden conjecture on permanents, a surprise at the time but signposted 20 years earlier by the proof by Marvin Marcus and Morris Newman of a weaker result. Egorychev and Falikman both used a simultaneity from geometry in their proofs! This was the Alexandrov-Fenchel inequality of 1936/37. Probably Thomas Clausen and Karl von-Staudt's simultaneous 1840 discovery of their eponymous theorem belongs also to the `unremarkably inevitable' category, given the ubiquity of the Bernoulli numbers in number theory following Euler, although a little detective work might be done: Clausen's contribution is in the form of a 10-line communiqu\'e appearing in Heinrich Christian Schumacher's journal {\em Astronomische Nachrichten}, signed `S' (presumably Schumacher), so it is unclear whether he had a proof. The theorem is also a 1911 repeat discovery, appearing without proof in Ramanujan's first published paper.

Naturally all mathematical discoveries, simultaneous or not, have a story to tell. The puzzle is to decide if there is a useful way of categorising these stories, or of extrapolating from them, or of accounting for them! Perhaps we can hypothesise a few.

Firstly we can identify cases where a simultaneous discovery is stimulated by an `unlocking idea'. Riemann's 1859 paper would fall into this category, as would Marcus and Newman's result on permanents. The exciting independent 2013 advances on small primes gaps of Zhang and Maynard where unlocked by a breakthrough of Goldston, Pintz and Y{\i}ld{\i}r{\i}m which came in 2005 after an almost twenty year stalemate.
Around 2003, no fewer than three independent proofs appeared that Moufang loops obey the Lagrange property: conversely to Cauchy's theorem, the order of a subgroup (or subloop) must divide that of its parent group/loop. All these proofs depended on a 70-page 1987 classification theorem of Peter Kleidman.

Unlocking can work in reverse, Michael O'Nan and Leonard Scott independently proved their classification theorem (the O'Nan-Scott Theorem) in 1979 as a theorem about maximal subgroups of the symmetric group. It was quickly harnessed as way of applying to primitive permutation groups the (still incomplete) Classification of the Finite Simple Groups (CFSG). Cameron recounts\cite{Cameron2} how ``Peter Neumann showed me that most of the [O'Nan-Scott] theorem is contained in Camille Jordan's {\em Trait\'e des Substitutions} from 1870."  Jordan's work had been overlooked: perhaps, Cameron suggests, ``the most plausible hypothesis is that, until we had CFSG, the theorem was not of much use, and was not considered worth pursuing."

CFSG is more than an `idea' of course - it counts as `big machinery'. In general it might be proposed that big machinery would act as a deterrent to simultaneity. CFSG, via the O'Nan-Scott Theorem, led to many conjectures being resolved: ``low-hanging fruit waiting to be picked", Cameron calls it\cite{Cameron2}; the circle of those able to do the picking would initially be small and would coordinate their investigations. If this sounds like a Mertonian reward system reinforcing itself then it should be noted that within two years of O'Nan and Scott's result, Cameron had written an `entry-level' survey paper\cite{Cameron1} making the low-hanging fruit accessible to the profession at large! Be that as it may, there are major results in permutation group theory, proved in the early 1980s by groups collaborating to exploit O'Nan-Scott, which remain reliant on CSFG to this day.

In graph theory the first big machinery came in the form of what Lov\'asz\cite{Lovasz} has called `the decomposition paradigm', of which the major example is the Robertson--Seymour Graph Minors Theorem. As with CFSG, the machinery was made accessible quite quickly\cite{Diestel} but it seems fair to say that those extending and using the paradigm  to solve deep problems remain a small and well-coordinated group. Here is an example: a well-known 1977 conjecture of Fiorini and Wilson asserts that every uniquely edge-3-colorable cubic planar graph on at least four
vertices contains a triangle. This conjecture is actually a simultaneity because S. Fisk made the same conjecture at the same time\cite{Fisk} (see Problem 11 in section I.6) but this is unsurprising since graph colourings were intensively studied at the time (the Four Colour Theorem was proved in 1976) and anyone at the cutting edge of the subject might have made this conjecture. What is interesting is that a proof came in 1998 using the decomposition paradigm, constituting the PhD of Tom Fowler, a student of Robin Thomas who is a leading exponent of the paradigm. The proof is long and difficult; it was announced by Thomas\cite{Thomas} in an overview paper but does not exist outside of Fowler's dissertation. Although the proof is computer-assisted nobody doubts its correctness. But equally, nobody is likely to look for a new proof outside of the decomposition paradigm, and the chances of a simultaneous proof in 1998 must have been very remote.

Despite a potential trend towards teamwork, data spanning the six decades up to the millenium\cite{Grossman} do not suggest a dramatic recent increase in authors per paper. But
beyond `big machinery' there is on the horizon an even bigger `threat' to simultaneity, the computerised proof. This is foreshadowed by the computer-assisted proofs of Fowler, and of Appel--Haken. It has come closer with recent formal (as opposed to computer-assisted) proofs of the Jordan Curve Theorem,\cite{Hales2} the Four-Colour Theorem,\cite{Gonthier1} the Feit--Thompson Theorem of finite group theory,\cite{Gonthier2} the Prime Number Theorem (both the elementary\cite{Avigad} and non-elementary\cite{Harrison} proofs) and Kepler's Conjecture\cite{Hales3}. Not many professional mathematicians yet take seriously the vision of Gowers, put forward in 2000,\cite{Gowers} of a time when computers both discover and proof all interesting theorems. But even a mathematical world in which, to be taken seriously, a theorem's proof had to be certificated as formally checked, would be a world in which mathematical activity was more regulated, more homogenised, more expensive. Theorems would be discovered and proved, even more than in the `big machinery' scenario, in large, coordinated groups (the publication describing the proof of Feit--Thompson, for instance, lists no fewer than fifteen co-authors!) The possibility of simultaneous discovery, and still less of repeated discovery, would surely be sadly diminished.

Having the same effect (but more congenial to the traditionally-minded mathematician) is the {\em Polymath} research initiative in which the Internet is used to focus the skills of hundreds of mathematicians on to a nominated open problem. The {\em Polymath 8} project, for example, resulted the publication,\cite{Polymath} under the name D.H.J. Polymath, of significant improvements to  bounds arising from Zhang and Maynard's small prime gaps breakthrough, improvements which could hardly have been achieved by one person working on their own, let alone two such simultaneously.

Mathematicians are curious about how they go about their work, even if they feel, with G.H. Hardy, that it is an idle curiosity. Simultaneous and repeated discovery, Merton's `multiples', are worth being curious about but not, for mathematicians, in the spirit of sociology, nor even of social history. It is more a matter of the solemn observation that `great minds think alike' leading to the question `why?'; or the observation that `it isn't worth re-inventing the wheel' leading to the question `why not'? Such investigations can reveal interesting information at the level of clever individuals proving clever theorems, even if the professional historian might more admire Stedall's ``... how, where, and why...".

These investigations may prove historical in a different sense, however: simultaneity and rediscovery might be consigned to history if the clever theorems are proved not by individuals but by clever teams backed by clever computer programmers. We might see simultaneous discovery replaced by races between teams competing to be the first to discover. At this point Merton's reward system view of science would indeed apply to mathematics.


\end{document}